\documentstyle[11pt]{article}
\renewcommand{\baselinestretch}{1.3}

\newtheorem {th}{Theorem}

\def\Cox{\hfill \Box}

\def\ee{\epsilon}
\def\E{{\bf{E}}}
\def\P{{\bf{P}}}
\def\N{\hbox{I\kern-.2em\hbox{N}}}
\def\R{\hbox{I\kern-.2em\hbox{R}}}

\def\|{\, | \, }

\def\0{\hat{0}}
\def\1{\hat{1}}

\begin{document}
\begin{titlepage}
\begin{center}
{\large \bf Sharpness of second moment criteria for branching and
tree-indexed processes} \\
\end{center}
\vspace{5ex}
\begin{flushright}
Robin Pemantle \footnote{Department of Mathematics, University of 
Wisconsin-Madison, Van Vleck Hall, 480 Lincoln Drive, Madison, WI 53706}$^,$
\footnote{Research supported in part by National Science Foundation Grant 
\# DMS 9300191, by a Sloan Foundation Fellowship and by a Presidential
Faculty Fellowship} ~\\
\end{flushright}

\vfill

\noindent{\bf ABSTRACT:}

A class of branching processes in varying environments is exhibited
which become extinct almost surely even though the means $M_n$ grow fast
enough so that $\sum M_n^{-1}$ is finite.  In fact, such a process
is constructed for every offspring distribution of infinite variance,
and this establishes the converse of a previously known fact: that
if a distribution has finite variance then $\sum M_n^{-1} = \infty$ 
is equivalent to almost sure extinction.  This has as an immediate
consequence the converse to a theorem on equipolarity of Galton-Watson 
trees.


\end{titlepage}

This note provides a class of examples of branching processes in  
varying environments (BPVE's) that die out almost surely even though 
the means grow relatively fast.  
\begin{th} \label{th 1}
Let $f$ be any offspring generating function with $f'(1) = m > 1$
and $f'' (1) = \infty$.  Then there is a sequence of positive
numbers $p_n \leq 1$ such that
\begin{eqnarray*}
(i) & ~ & \sum_{n=1}^\infty \prod_{k=1}^n (m \cdot p_k)^{-1} < \infty
   ~~\mbox{  and} \\[2ex]
(ii) & ~ & \mbox{The BPVE with offspring generating functions }
   f_n (z) = f(1 - p_n + p_n z) \\
&& \mbox{dies out almost surely.}
\end{eqnarray*}
\end{th}

The expected size of the $n^{th}$ generation of the BPVE in~$(ii)$ 
is given by $$\E Z_n = \prod_{k=1}^n (m \cdot p_k) $$
which is why condition~$(i)$ is a growth condition on the means.   
The proof of this theorem is easy, and the exposition will focus mostly 
on saying why the theorem is interesting.  There are two reasons, one
having to do with branching processes and one to do with tree-indexed
processes.  The motivation coming from branching processes is more
straightforward.  

A branching process in a varying environment (BPVE) is defined
by a sequence of offspring generating functions 
$$ f_n (z) = \sum_{k=0}^\infty q_{n,k} z^k $$
where for each $n$, the nonnegative real numbers $\{ q_{n,k} \}$
sum to 1.  From the sequence $\{ f_n \}$ a random tree $\Gamma$
is constructed as follows.  The root has a random number $Z_1$ of
children, where $\P (Z_1 = k) = q_{1,k}$.  Each of these first-generation
individuals has a random number of children, these random numbers
$X_1 , \ldots , X_{Z_1}$ being IID given $Z_1$ and satisfying
$\P (X_1 = k) = q_{2,k}$.  This continues in the same manner,
so that if $Z_n$ is the total number of individuals in generation $n$,
then the numbers of children of each of these $Z_n$ individuals are
IID, being equal to $k$ with probability $q_{n+1,k}$.  

The mean number of children of an individual in generation $n$ is
$f_n'(1)$ and therefore the expected number of individuals
in generation $n$ is given by
$$M_n = \prod_{k=1}^n f_k' (1) .$$
A {\em Galton-Watson} or simple branching process is one where
the environment does not vary, i.e., $f_n = f$ for
all $n$.  In this case $f' (1) \leq 1$ is necessary and sufficient 
for almost sure extinction:
$$f'(1) \leq 1 \Longleftrightarrow \P (Z_n \rightarrow 0) = 1 .$$
In the varying case, it is possible for $M_n$ to grow without
growing exponentially, and one may ask whether the growth 
rate of $M_n$ determines whether $\P (Z_n \rightarrow 0)$ is 
equal to 1.  Under the assumption that $\liminf_{n \rightarrow \infty}
1 - f_n(0) - f_n'(0) > 0$ (a weak nondegeneracy condition saying the
probability of at least two children is bounded away from zero),
Agresti (1975) shows that
$$\sum M_n^{-1} = \infty ~~\mbox{ implies almost sure extinction } . $$
On the other hand, a second moment condition is needed for
the converse.  The second moment of the number of children in
generation $n$ is $f_n'' (1) + f_n' (1)$.  Agresti shows
that if $\sup_n f_n'' (1) < \infty$ then 
\begin{equation} \label{eq 1}
\sum M_n^{-1} < \infty ~~\mbox{ implies a positive probability
   of non-extinction}. 
\end{equation}
See also Theorem~4.14 of Lyons (1992).  

A natural class of BPVE's are those obtained from a single offspring
generating function $f$ by killing individuals in generation $n$
independently with probabilities $1 - p_n$.  Think of this as
modeling a genealogy where the branching mechanism remains the
same from generation to generation but the hospitality of the 
environment varies.  Formally, $f_n (z) = f(1 - p_n + p_n z)$, so 
the expected generation sizes are $M_n = \prod_{k=1}^n (m \cdot p_k)$.  
The second moment condition $\sup_n f_n'' (1) < \infty$ 
is equivalent in this case to $f''(1) < \infty$.  Theorem~1 shows
that this condition is necessary as well as sufficient 
for~(\ref{eq 1}) to hold: for any $f$ with infinite second moment,
some BPVE of the form $f_n (z) = f (1 - p_n + p_n z)$
becomes extinct almost surely even though $\sum M_n^{-1} < \infty$.

Theorem~1 may also be viewed as a fact about tree-indexed processes.
A process indexed by a tree $\Gamma$ is simply a set of IID real 
random variables $\{ X(v) \}$ indexed by the vertices of $\Gamma$.  
Let $B \subseteq \R^\infty$ be closed in the product topology.
The following notion of polar sets for tree-indexed processes was
first defined by Evans (1992).  \\[3ex]
{\bf Definition}: {\em The set $B$ is polar for $\Gamma$ (and
for the common distribution of the variables $\{ X(v) \}$)
if and only if the probability is zero that there exists an 
infinite self-avoiding path $v_0 , v_1 , v_2 , \ldots$ from 
the root of $\Gamma$ satisfying $(X(v_1) , X(v_2) , \ldots ) 
\in B$.} \\[3ex]
Trees with the same polar sets are denoted {\em equipolar} by 
Pemantle and Peres (1994).  In particular,
letting $\{ X(v) \}$ be uniform on the unit interval and
letting $B = \{ (x_1 , x_2 , \ldots ) : \forall n \; x_n \leq p_n \}$,
one sees that equipolar trees $\Gamma_1$ and $\Gamma_2$ are 
percolation equivalent, meaning that:
\begin{quote}
If vertices of both trees are removed independently with the
survival probability $p_n$ of a vertex being the same for all
vertices in generation $n$ of either tree, then the root of $\Gamma_1$
has positive probability of being in an infinite component of
surviving vertices if and only if the root of $\Gamma_2$ has positive
probability of being in an infinite component of surviving vertices. 
\end{quote}
There is not space here for a substantial discussion 
of equipolarity, but the reader is referred to Pemantle and 
Peres (1994), wherein it is shown that equipolar trees behave similarly 
for a variety of common probability models (other than percolation) 
involving trees, including maximal displacements of branching random 
walks, nonextinction probabilities for branching random walks with
absorption, survival of certain BPVE's, growth rates of 
first-passage percolation clusters, and capacities of fractal 
sets in Euclidean space defined by interpreting the tree as a 
base-$b$ expansion of a closed subset of the unit cube.  
Equipolarity results for random trees lead to Peres' (1994) derivation 
of Fitzsimmons and Salisbury's (1989) capacity criteria for multiple 
points of Brownian motions and to a general capacity-theoretic 
framework for intersection properties of random sets.  Given
that equipolarity is a useful notion, the effort to understand
which trees are equipolar should seem justified.  

In Pemantle and Peres (1994) it is shown that trees which arise from 
Galton-Watson processes with respective offspring generating functions 
$f$ and $g$ are almost surely equipolar provided they have the
same mean growth $f'(1) = g'(1)$,  and that each has a finite
variance: $f'' (1) < \infty$ and $g'' (1) < \infty$.  It is
also shown in the preprint version that the second moment assumption is
almost sharp in the sense that if $f'' (1) < \infty$ but
the distribution defined by $g$ fails to have a $2 - \ee$ 
moment for some positive $\ee$, then there is some set $B$ 
which is almost surely polar for a Galton-Watson tree with
offspring generating function $g$ but almost surely nonpolar 
for a Galton-Watson tree with offspring generating function $f$.
Theorem~1 improves this to a sharp result, namely that 
whenever $g$ fails to have a second moment, the set 
$$B = \{ (x_1 , x_2 , \ldots ) : \forall n \; x_n \leq p_n \} $$
defined from the sequence $\{ p_n \}$ in the conclusion of the 
theorem is polar for almost every Galton-Watson tree with offspring
generating function $g$.  (The fact that it is nonpolar for 
almost every Galton-Watson tree with offspring generating function
$f$ follows from $\sum M_n^{-1} = \infty$.)

Having given motivation for the theorem, I now give the proof, 
which is based on a well-known result of Kesten, Ney and Spitzer
(proved with a third moment assumption by Kolmogorov).
\begin{th}[K-K-N-S] \label{th K}
Suppose $g$ is an offspring generating function for a 
critical Galton-Watson process, i.e., $g' (1) = 1$.  Let
$\sigma^2 = Var (Z_1) = g''(1) \leq \infty$.  Then
$$\lim_{n \rightarrow \infty} n \P (Z_n > 0) = {2 \over \sigma^2} .$$
\end{th}

\noindent{\sc Proof:}  See Kesten, Ney and Spitzer (1966) or Lyons,
Pemantle and Peres (1994).   $\Cox$

\noindent{\sc Proof of Theorem~1:}  Fix an offspring generating
function $f$ with $m := f'(1) < f''(1) = \infty$.  Let $g (z) = 
f(1 - 1/m + z/m)$ so that $g$ is an offspring generating function 
satisfying the hypotheses of K-K-N-S theorem.  Applying the
theorem to the probabilities $P_n$ of the critical branching process
with offspring generating function $g$ surviving $n$ levels, we
see that we may choose for every $n$ an $L_n$ such that for all $k \geq L_n$,
$$P_k < 4^{-n} k^{-1} .$$
Define sequences $\{ t_n \}$ and $\{ u_n \}$ recursively as follows.
Let $u_0 = 0$ and let $t_1$ be the least positive integer for which 
$m^{t_1} > L_1$.  For each $n \geq 1$, let 
$$K_n = \prod_{j=1}^n m^{t_j - u_{j-1}} ,$$
let
$$u_n = t_n + \lceil 2^{-n} K_n \rceil $$
and let $t_{n+1}$ be the least integer $k > u_n$ for which 
$$m^{k - u_n} \prod_{j=1}^n m^{t_j - u_{j-1}} > 2^{n+1} L_{n+1} .$$
This ensures that $K_{n+1} > 2^{n+1} L_{n+1}$.
Let $p_n = 1/m$ whenever $t_j \leq n < u_j$ for some $j$, and
let $p_n = 1$ whenever $u_j \leq n < t_{j+1}$ for some $j$.
When $t_j \leq n \leq u_j$, the $n^{th}$ generation has expected size
$\E Z_n = \prod_{k=1}^n (m \cdot p_k) = K_n$.    

To verify condition~$(i)$ of Theorem~1, observe first that
if $b_n$ is the $n^{th}$ positive integer $j$ (in ascending order)
for which $p_j = 1$, then 
$$\sum_{n=1}^\infty \prod_{k=1}^{b_n} (m \cdot p_{k})^{-1}
   = \sum_{n=1}^\infty m^{-n} < \infty .$$
Thus it suffices to show that
\begin{equation} \label{eq 11}
\sum_{n : p_n = {1 \over m}} \;  \prod_{k=1}^n (m \cdot p_k)^{-1} < 
   \infty . 
\end{equation}
Write this sum as
$$\sum_{n=1}^\infty \sum_{j=t_n}^{u_n - 1} \prod_{i=1}^j 
      (m \cdot p_i)^{-1}  
   = \sum_{n=1}^\infty (u_n - t_n) K_n^{-1} .$$
By construction, $u_n - t_n < 2^{-n} K_n + 1$, and since
$K_n^{-1} < m^{-n}$, the sum in~(\ref{eq 11}) is finite.

To verify condition~$(ii)$, first write
\begin{equation} \label{eq 2}
\P (Z_{u_n} > 0) = \E (\P (Z_{u_n} > 0) \| Z_{t_n}) \leq
   r_{u_n - t_n} \E Z_{t_n} 
\end{equation}
where $r_n$ is the probability that the critical branching
process with offspring generating function $g$
survives to the $n^{th}$ generation.  Since $u_n - t_n \geq
2^{-n} K_n > L_n$, the definition of $L_n$ then forces
$$r_{u_n - t_n} < 4^{-n} (u_n - t_n)^{-1} .$$
Plugging this into~(\ref{eq 2}) and using $\E Z_{t_n} = K_n$
then gives $\P (Z_{u_n} > 0) < 4^{-n} [2^{-n} K_n]^{-1} K_n$,
and hence
$$\lim_{n \rightarrow \infty} \P (Z_{u_n} > 0) = 0 .$$
Thus the BPVE dies out almost surely.   $\Cox$

\renewcommand{\baselinestretch}{1.0}\large 
\vspace{.4in}

\today

\end{document}